\documentstyle[twoside,amssymb,12pt]{article}
\title{\sc Regular minimal nets on  surfaces 
 of constant negative curvature
}

\author{{\sc Alina Vdovina} \\
\small  {Khmelnitskiy Thechnical University }  \\
\small  { e-mail: \ \ alina@alpha.podol.khmelnitskiy.ua}
\and 
{\sc E.N. Selivanova}\thanks{Supported by DAAD.}\\
\small  { Department of Geometry, Nizhny Novgorod State 
Pedagogical University}\\
\small  { 603000 Russia, Nizhny Novgorod, ul. Ulyanova 1}\\
\small  {e-mail: \ \ lena@moebius.mathematik.uni-tuebingen.de} 
}

\begin{document}
\date{}

\maketitle

\noindent
{\bf Abstract.} {\small The problem of classification of closed local minimal 
nets on surfaces of constant negative curvature has been formulated in \cite{Fo1}, \cite{Fo2} 
in the context
of the famous Plateau problem in the one-dimensional case. In \cite{SeV} an asymptotic for 
$\log \sharp (W^r(g))$ as $g\to +\infty$ where $g$ is  genus and $W^r(g)$ is the set 
of regular single-face closed local minimal 
nets on surfaces of curvature $-1$ has been obtained. It has been shown that $\sharp (W^r(g))$
is equal to the number of classes of topological equivalence of all  single-face
closed local minimal nets on surfaces of curvature $-1$.
In this paper we prove an  asymptotic for $\sharp (W^r(g))$ as $g\to +\infty$ and 
construct some examples of $f$-face, $f>1$  nets.}

\thispagestyle{empty}

\section{\bf Introduction}

Let $(M^2, ds^2)$ be a Riemann surface of constant curvature. 
{\it Closed local minimal net} $\Gamma$ on $M^2$ is an embedded 
cubic graph consisting of
segments of geodesics where the angles at any vertex are equal to 
$2\pi/3$.
A local minimal net is {\it regular} if all segments of geodesics are 
equal.

Two local minimal nets are {\it toplogically equivalent}
if there is a homeomorphism of surfaces 
preserving the orientation and transferring nets one into another.

Two local minimal nets are {\it isomorphic}
if there is an isometry of surfaces  transferring these nets one into 
another.

Let $\Gamma$ be a net on $M^2$. 
If the difference $M^2\setminus \Gamma$ is 
homeomorphic to one disk then
the net $\Gamma$ is {\it single-face}.
If the difference $M_g^2\setminus \Gamma$ is homeomorphic to $f$  disks 
then
 net $\Gamma$ is  $f$-{\it face}.

The problem of classification of closed local minimal 
nets on surfaces of constant negative curvature has been formulated in \cite{Fo1} in the context
of the famous Plateau problem in the one-dimensional case. 

Full description of such nets in the case of curavature $1$ one can find in \cite{He}.
There are many interesting applications of this result, see \cite{AT}, \cite{Fo1}, \cite{Fo2}.

We consider the case of orientable surfaces of curvature $-1$. 
In \cite{SeV} an asymptotic for 
$\log \sharp (W^r(g))$ as $g\to +\infty$ where $g$ is  genus and $W^r(g)$ is the set 
of regular single-face closed local minimal 
nets on surfaces of curvature $-1$ has been obtained. It has been shown that $\sharp (W^r(g))$
is equal to the number of classes of topological equivalence of all single-face
closed local minimal nets on surfaces of curvature $-1$.
In this paper we prove an  asymptotic for $\sharp (W^r(g))$ as $g\to +\infty$ and 
show that a description 
of all single-face closed  minimal nets on surfaces of curvature $-1$ can be obtained d
ue to a description 
of the set 
of regular single-face closed  minimal 
nets in terms of polygons on the hyperbolic plane. Then we
construct some examples of $f$-face nets, $f>1$.

\section{\bf Single-face nets}
\bigbreak
\noindent
In order to describe the set of the single-face regular local minimal nets on $M_g^2$
 will use the following algebraic construction.

Let us call a  word $U$ of length $12g-6$ a Wicks form of genus $g$ if it satisfies
 the following conditions

1. each letter has been used twice as  $x$ and $x^{-1}$,

2. there is no pair of distinct,
noninverse letters $x,y$ which appear in $U$ only in subwords 
$(xy)^{\pm 1}$.

If letters $x,y$  appear in $U$ only in subwords 
$(xy)^{\pm 1}$, we can replace $(xy)^{\pm 1}$ by $z^{\pm 1}$.
This operation we will call {\it reduction.}

Two words $U_1$ and $U_2$ are {\it equivalent}, if  
$U_1$ can be obtained from $U_2$  by permutation of 
letters,
by substitution of some of them for inverse ones and by cyclic
permutation.

It turns out that there is an {\it inductive
method} to construct all  words $U(g+1)$ from $U(g)$.
We can describe the {\it transformations} by which  all words 
$U(g+1)$
can be obtained from  words $U(g)$.

Let us rewrite  $U$ changing $a$ into  $a_1a_2$ 
  and put  the fragment
$$id^{-1}ef^{-1}dhe^{-1}fh^{-1}i^{-1}$$ between
$a_1$ and $a_2$ or  between $a_2^{-1}$ and $a_1^{-1}$. This 
transformation
is called  $\alpha$-transformation.

Another transmormation is to  rewrite   $U$ changing $a$ into 
 $a_1a_2$ and $b$ into  $b_1b_2$
 and put  the fragment   $$d^{-1}ef^{-1}d$$ between
$a_1$ and $a_2$ or  $a_2^{-1}$ and $a_1^{-1}$ and put the fragment
$$he^{-1}fh^{-1}$$ between
$b_1$ and $b_2$ or  $b_2^{-1}$ and $b_1^{-1}$.This transformation
is called by $\beta$-transformation.

The third transformation is the following. We 
rewrite $U$ changing $a$ into $a_1a_2$,  $b$ into $b_1b_2$  and
$c$ into $c_1c_2$ i.e. $$U^{\prime}=U_1a_i^\epsilon 
a_j^\epsilon U_2b_k^\delta b_l^\delta U_3c_m^\psi c_p^\psi U_4$$ (an 
expression of type 
$a_i^\epsilon a_j^\epsilon$
 means   the  $a_1a_2$ or $a_2^{-1}a_1^{-1}$,
where  $\epsilon=1$, $i=1$, $j=2$ or $\epsilon=-1$, $i=2$, $j=1$)
 We'll say that  the word  $$V=U_1a_i^\epsilon d^{-1}eb_l^\delta
U_3c_m^\psi f^{-1}d
 a_j^\epsilon U_2b_k^\delta e^{-1}fc_p^\psi U_4$$ is obtained from 
$U$
by  $\gamma$-transformation.

\medbreak
\noindent
{\bf Lemma, \cite{6}.} Word $U(g)$ corresponds to an embedded cubic graph 
with $4g-2$ verties
on a $M_g^2$ if and only if this word  can be obtained from
$$abca^{-1}b^{-1}c^{-1}$$   by a sequence of   $(g-1)$
$\alpha$, $\beta$ or $\gamma$-transformations.
\bigbreak
\noindent
{\bf Theorem 1.1} For $\sharp W^r(g) $ the following asymptotic 
$$\sharp  W^r(g) =\frac{(6g-4)!}{(3g-2)!g!12^g}(1+o(1))$$ holds.
\medbreak
\noindent
{\it Proof.} As mentioned above, $\sharp  W^r(g)$ equals  the number of classes of topological
 equivalence of the single-face nets on $M_g$. First we shall show that the problem of topological 
classification
of  single-face local minimal nets on $M_g^2, g>1$ is equivalent to the problem of
topological classification of embeddings of cubic graphs to a
surface $M_g^2$.

Let us consider a single-face  embedded cubic graph $P$ in a surface 
$M_g^2$, 
label the edges of $P$ by different
 letters $\{ {a,b,...} \}$ of an alphabet $A$.
Let us choose also an orientation of  any edge of $P$.
The difference $M_g^2\setminus P$ consists of a disk with 
a word $U$ on the boundary of this disk. The number
of the different letters (or the number of edges of the graph)
is $6g-3$. 
In the word $U$ each letter has been used twice.
Taking into account the
orientation of the edges we'll write $x$ and $x^{-1}$ if 
the orientations of the
corresponding segments are different. There is no pair of distinct,
noninverse letters $x,y$ which appear in $U$ only in subwords 
$(xy)^{\pm 1}$. 
   Our classification theorem
is based on the description of words $U$ from Lemma befor the 
theorem.

There is a realization of local minimal nets for any 
equivalence
 class  of cubic graphs on the $M_g^2$ as a nontrivial subset
of this class. (One may consider a regular polygon on the hyporbolic plane labeed according the corresponding word $U$).

Thus, we have now to describe the set of equivalence classes of   
graphs.
This means to describe all possible words $U$ up to equivalence, see 
above.

Let us consider an oriented connected cubic graph $\Gamma$ with $2k$ 
vertices
and $3k$ edges labeled by different 
letters of an alphabet $A$ with specified cyclic orderings of the 
edges at any
vertex. At any vertex $v$ of graph we will write a cyclic ordered 
triples
of letters so that a letter $x$ has exponent $1$ if the edge $x$ 
enters
in   $v$ and $-1$ if $x$ leaves $v$. An {\it ordered graph} is a 
graph 
with specified cyclic oderings (the triples) at the edges at each 
vertex.
We will call the set of elements $x, x^{-1}, x\in A$ by $2Q$. 
 The vertices of a graph may be considered  as  sets into which $2Q$
is partitioned.
Let us consider a permutation 
$P$ on $2Q$ such that the cyclic order at any vertex 
determines a cycle in the  permutation $P$ and an inverse 
 $\xi$ transferring $x$ to $x^{-1}$ for any $x$ from $2Q$.
If  permutation $P\xi$ has $f$ cycles then  permutation $P$
determines an embedding of graph $\Gamma$ in the sphere with $g$
handles where there are  $f$ faces, $2k$ vertices and $3k$ edges,
satisfying the equality $k=2g+f-2$, \cite{7}, i.e. there is an 
$f$-face embedding of genus $g$. (Note that for $f=1$
permutation $P\xi$ is the considered above Wicks form $U$). Two
ordered graphs $(Q,P)$ and $(Q',P')$ are {\it isomorph} if there is 
a one-to-one  mapping $\Phi$ such that 
$$\Phi[P(\beta)]=P'[\Phi(\beta)], \Phi(\beta^{-1})=\Phi(\beta)^{-1}$$
for any $\beta$
from $2Q$. Thus the {\it  automorphism} group of $P$ (denoted by $aut 
(P)$)
is the set of permutations on $2Q$ which commute with $P$ and $\xi$.

A connected ordered graph is rooted
if one element $\beta$ of $2Q$ is distinguished. A {\it root 
isomorphism}
from one ordered graph $(E,P, \beta)$ onto another one $(E',P', 
\beta')$
is an isomorphism $\Phi$ from $(E,P)$ onto $(E',P')$ such that 
$\Phi(\beta)=\beta'$.
Thus a {\it root automorphism } on $(E,P,\beta)$ is a permutation 
on $2E$ which commutes with $P$ and $\xi$ and preserves $\beta$.
(Corresponding Wicks form we will call a rooted Wicks form.)
In \cite{7} has been shown that {\it the number of essentially different 
ways to root 
the ordered graph} $(E,P)$  that is, the number of root-isomorphism 
classes of rooted ordered graphs which are isomorphic as ordered 
graphs to
$(E,P)$ is $6k/|aut (P)|$. Then 
$$1\le t\le 6k. \eqno (1)$$
Let $R(1,g)$ be the number of rooted ordered graphs on $4g-2$ 
vertices,
$6g-3$ edges such permutation $P\xi$ consists of one cycle. 

The following formula has been obtained in \cite{7}:
$$R(g)=
 {(12g-6)(6g-4)!\over{(3g-2)!g!12^g}}.\eqno (2)$$
Let $N(1,g)$ be the number of rooted ordered graphs on $4g-2$ 
vertices,
$6g-3$ edges such permutation $P\xi$ consists of one cycle. Then
$R(g)/t \leq N(1,g) \leq R(g)$.
There is a one-to-one correspondence between  single-face 
genus $g$ embeddings and ordered graphs where permutation $P\xi$ 
consists
of one $12g-6$. Then the number of genus $g$ embeddings equals to 
$N(1,g)$.
Taking into account (1) and (2) we obtain  
 $${(6g-4)!\over{(3g-2)!g!12^g}} \leq N(1,g) \leq$$
 $$\leq  {(12g-6)(6g-4)!\over{(3g-2)!g!12^g}}.\eqno(3)$$
Thus $$N(1,g)\ge \frac{R(g)}{12g-6}.$$
Let us prove $N(1,g)\le \frac{R(g)}{12g-6}(1+o(1))$.

\bigbreak
Denote $N_1(g)$ the number of the Wicks forms with trivial 
automorphisms
and $N_k(g)$ the number of the Wicks forms with automorphisms
of order $k$.
It is obvious, that $$N_1(g)\le R(g)/12g-6.$$
So, we have to prove, that 
\begin{equation}
\sum\limits^{12g-6}_{k=2}N_k(g)=
o(R(g)/12g-6).
\label{sss}
\end{equation}
Let $x$ and $y$ be some letters of $U$ where $U=U_1xU_2y$.
We define the distance between $x$ and $y$ (it is noted as $d(x,y)$) 
for
letters $x$, $y$ of $U$ to be $\min (|U_1|, |U_2|)$.
We will estimate the number of boundary walks of single-face
embeddings of cubic graphs on $4g-2$ vertices.
It is obvious, that a boundary walk of a single-face embedding
of a cubic connected graph on $4g-2$ vertices can be written
as a word $U$ with length $12g-6$.  The word $U$  
will be a genus $g$ Wicks form.
Let us note, that if edges $a,b,c$ are leaving a vertex $v$ of $G$,
then there are subwords $a^{-1}b,b^{-1}c$ and $c^{-1}a$ in $U$.
\medbreak
\noindent
{\bf Lemma 1.} One may estimate the number of the Wicks forms
 with 
automorphisms of the second order by the following expression 
$$N_2(g)\le (g+1) \prod_{m\in N, 12g-6-7m>0} (12g-6-7m).$$ 
\medbreak
\noindent
{\it Proof.}  We will consider a $12g-6$-gon
and try to construct  a Wicks form $U$ with an automorphism
$\xi$ of the second order.
 It is known, (see,for example \cite{6}), that every
multiple edge $(b,c)$ gives subwords of $abca^{-1}$ and 
$db^{-1}c^{-1}d^{-1}$
kind  in $U$ or subword of $eabca^{-1}db^{-1}c^{-1}d^{-1}e^{-1}$ 
kind.
In the begining of the proof we consider the case, when $U$
has $l$ multiple edges such that $db^{-1}c^{-1}d^{-1}= 
\xi(abca^{-1})$.
We will say, that the multiple edge $(b,c)$ is {\it diagonal multiple
edge.}

On the first step we fix  words
   $a_1b_1c_1a_1^{-1}$ and $d_1b_1^{-1}c_1^{-1}d_1^{-1}$
corrisponding to a multiple edge $(b_1,c_1)$ ,
$d_1b_1^{-1}c_1^{-1}d_1^{-1}= \xi(a_1b_1c_1a_1^{-1})$. So,
$8$ letters are fixed before the second step and there are no more
than $12g-6-8$ possibilities for the words corresponding to the 
second
multiple edge $(a_2,b_2)$. There are 16 letters are fixed before the 
third step
and there are no more than $12g-6-16$ possibilities for the words 
corresponding to the third multiple edge. The procedure we will call
$\alpha$- procedure and $8k$ letters have been already fixed after 
$k$-step of
the $\alpha$- procedure. The $\alpha$- procedure has l steps and
there are $12g-6-8l$ "free" places after it.
Let's describe a $\beta$-procedure. Let's note, that if letter
$x$ is fixed, then $x^{-1}$, $\xi(x)$ and $\xi(x^{-1})$ are fixed 
too.
Let $z$ be a letter, which have been already fixed, but there is
a free place near $z$. We continue $l$-step of $\alpha$- procedure
by fixing elements $x_1$ and $y_1$ such that $z^{-1}x_1$ and 
$y_1^{-1}z$ are
subwords of $U$ and $x_1$ and $y_1$ have not been fixed by $\alpha$- 
procedure.
So, there are no more, than $12g-6-8l-2$ possibilities for
the subword $x_1^{-1}y_1$ in $U$. The places for
$\xi(x_{1}), \xi(y_{1}), 
\xi(x_{1}^{-1}),\xi(y_{1}^{-1})$ are determined, because $U$
has an automorphism of the second order. Let $w$ be a letter,
which has already fixed, but there are free places after $w$
and before $w^{-1}$. We will put $x_2$ after $w$ and $y_2$
before $w^{-1}$. We define $k$-step of the $\beta$-procedure
so that elements
$x_k, y_k, \xi(x_{k-1}), \xi(y_{k-1}), x_{k-1}^{-1}y_{k-1}^{-1},
\xi(x_{k-1}^{-1}),\xi(y_{k-1}^{-1})$ are fixed on the $k$-step.
Let $A_k$ be the set of all letters of $U$ which have been already 
fixed 
after the $k$-step. So, $|A_k|=|A_{k-1}|+|N(k)|$, where $N(k)$ is a 
number
of letters, which can be fixed on the $k$-step.
By definition of $k$-step, the set $N(k)$ consists of the elements
$x_k, y_k, \xi(x_{k-1}), \xi(y_{k-1}), x_{k-1}^{-1}y_{k-1}^{-1},
\xi(x_{k-1}^{-1}),\xi(y_{k-1}^{-1})$. According our construction
these elements has been not fixed before $k$-step. It easy to see, 
that
$N(k) \geq 7$, if $x_{k-1}$ and $y_{k-1}$ are not multiple
edges with $x_{k-1}= \xi(x_{k-1})$,$y_{k-1}= \xi(y_{k-1})$.

We will continue the $\beta$-procedure while it is possible.
So, the general procedure of constructing $U$ consists of two 
inductive
procedures( $\alpha$ and $\beta$).  
We will say, that the procedure of constructing of a cubic with
$l$ diagonal multiple edges is {\it successful}, if we obtain a Wicks
form with exactly $l$ diagonal multiple edges after it. 
Since we have no more than $12g-6-7k$ possibilities of choice 
in the $k$ step of the general procedure and no less than seven
letters are fixed in the each step,
then
$$N^l \leq (12g-6-8) \prod_{m \geq 2 , 12g-6-7m>0} (12g-6-7m) ,
l=1, \ldots, g,$$
where $N^l$ is the number of cubic words with an automorphism 
of the second order
with $l$ diagonal multiple edges (it is easy to see, that
a cubic word cannot have more than g multiple edges), $l=0, 
\ldots,g$.

 If $U$ has no multiple edges, we can start $\beta$-procedure with
 a letter $c$ such that $d(c,c^{-1})=6g-4$.
Since $Aut(U)=2$  the letter $c$ there exists. In this case in the 
first
step we fix only six letters(in the $k$-step, where $k \geq 2$ we fix
no less than seven letters by definition of the $\beta$-procedure),
$$N^0 \leq (12g-6-6)\prod_{m \geq 2 , 12g-6-7m>0} (12g-6-7m). $$
So, $$N_2(g)=N^0+N^1+ \ldots +N^g \leq 
(g+1) \prod_{m\in N, 12g-6-7m>0} (12g-6-7m)$$
The Lemma 1 is proved.
\medbreak
\noindent
{\bf Lemma 2.} One may estimate the number of the Wicks forms with 
automorphisms of  order $3$ by the following expression 
$$N_3(g)\leq 2g^2 \prod_{m\in N, 12g-6-7m>0} (12g-6-7m).$$ 
\medbreak
\noindent
{\it Proof.} 
Let $U$ be a Wicks form  with automorphism $\xi$
of the third order. Let $v$ be a vertex and edges $a,b,c$ are leaving
it. We will say, that $v$ is automorphic, if $\xi(a^{-1}b)=b^{-1}c$
and $\xi(b^{-1}c)=c^{-1}a$ or $\xi(a^{-1}b)=c^{-1}a$ and
$\xi(c^{-1}a=b^{-1}c$.

Let $v_1$ be an automorphic vertex, edges $a_1,b_1,c_1$ are leaving 
it
and the edge $b_1$ enters a vertex $w_1$. It is obvious,
that $w_1$ cannot be automorphic. Edges $b_1,b_2,b_3$ enter the 
vertex
$w_1$ and $b_2$ leaves a vertex $v_2$. The vertex $v_2$ can be 
automorphic. If $v_2$ is automorphic, we will say, that $v_2$
is {\it neighbor} of $v_1$. So, two automorphic vertices $v$ and $w$
are {\it neighboring}, 
if the shortest path between $v$ and $w$ consists of two edges.

Vertices, which are ends of an edge $x$ we will denote
by $x_{-}$ and $x_{+}$, where $x$ leaves $x_{-}$ and enters $x_{+}$.
Let us prove, if $|U| \geq 30(g \geq3)$, then every 
automorphic vertex can have only
 one neighbor. Let edges $a,b,c$ be leaving a vertex $v$,
edges $d,f,e$ be leaving a vertex $w$ and $v$ and $w$ are 
neighboring automorphic vertices, then $a_{+},b_{+},c_{+}$
and $d_{+},f_{+},e_{+}$ are equal sets of vertices of $\Gamma$,
that is if $b_{+}=e_{+}$, then $a_{+}=d_{+},c_{+}=f_{+}$ or
$a_{+}=f_{+},c_{+}=e_{+}$. Without loss of generality we can
assume, that $a_{+}=d_{+}$, $b_{+}=e_{+}$ and $c_{+}=f_{+}$.
Since all vertices of the graph $\Gamma$ are degree three, there are
some edges $g,h,i$ such that $a_{+}=d_{+}=g_{+}$, $b_{+}=e_{+}=h_{+}$ 
and $c_{+}=f_{+}=i_{+}$. If $g_{+}$ is automorphic, then
$g_{+}=h_{+}=i_{+}$ and $\Gamma$ has only three vertices.
So, if $|U| \geq 30(g \geq3)$, then every 
automorphic vertex can have only one neighbor.

 Let $U=U_1U_2U_3$ be a Wicks form with
automorphism $\phi$ of the third order with $t$ single automorphic
vertices (without neighbors) and $s$ pairs of automorphic vertices,
such that $|U_1=|U_2|=|U_3|$ and $U_2= \phi(U_1), U_3= \phi^2(U_1)$.

Now we describe a transformation of $ \psi$-type. Let $U_1=x_1x_2 \ldots
x_k$. First, we delete all letters from all automorphic vertices.
Then, if both $x_i$ and $x_i^{-1}$ belong to $U_1$, then we leave
$x_i$ and $x_i^{-1}$ without changing. If $x_i$ belongs to $U_1$,
but $x_i^{-1}$ belongs to $U_2$ or $U_3$, we replace the image of $x_i^{-1}$
 in $U_1$ by $x_i^{-1}$. After all reductions and cancellations
we will obtain a word $W$, which is a Wicks form of genus $p$,
where $p \le (g+2)/3$. We will say, that $W$ is obtained from $U$
by a transformation of $ \psi$-type. On the other hand, 
there are no more, than $$3^{12p} {(12p-6)(12p-2) \ldots (12p-6+4s)
(12p-2+4s) \ldots (12p-6+4s+4t)\over{s!t!}}$$ possibilities to obtain
a genus $g$ Wicks form of the third with $t$ single automorphic
vertices and $s$ pairs of automorphic vertices from every rooted Wicks form
of genus $p$. The number of nonequivalent rooted cubic words
is $$ {(12p-6)(6p-4)!\over{(3p-2)!p!12^p}}.$$ So, the number of 
genus $g$ Wicks forms of the third with $t$ single automorphic
vertices and $s$ pairs of automorphic vertices , is no more, than
$$S(p)=3^{12p} (12p-6)(12p-2) \ldots (12p-6+4s)
(12p-2+4s) \ldots $$ $$(12p-6+4s+4t)(12p-6)(6p-4)!\over{s!t!(3p-2)!p!12^p}.$$
It is easy to see, that 
$$S(p) \le  2g^2 \prod_{m\in N, 12g-6-7m>0} (12g-6-7m).$$
\medbreak
\noindent
{\bf Lemma 3.} One may estimate the number of the Wicks forms with 
automorphisms of  order  $5$ by the following expression 
$$N_5(g)=o(R(g)/12g-6).$$ 
\medbreak
\noindent
{\it Proof.} It is easy to see, that a Wicks form with automorphism
of order five cannot have any "diagonal" edges and automorphic
vertices. So, every vertex $v_1$ has four different images.

If edges $a_1, b_1, c_1$ are leaving a vertex $v_1$, then
there are subwords $a_1^{-1}b_1, b_1^{-1}c_1$ and $c_1^{-1}a_1$
 in the Wicks form. Now we will describe an inductive procedure
of constructing a genus $g$ Wicks with automorphism
of order five. At the first step we fix subwords $a_1^{-1}b_1, 
b_1^{-1}c_1$ and $c_1^{-1}a_1$ and their images.
There are no more than $12g-6$ possibilities for fixing
of subword $a_1^{-1}b_1$, there are no more than $12g-6-2$
possibilities for fixing of subword $b_1^{-1}c_1$ and
there are no more than $12g-6-4$
possibilities for fixing of subword $c_1^{-1}a_1$.
The places for images of  $a_1^{-1}b_1, 
b_1^{-1}c_1$ and $c_1^{-1}a_1$ are determined and $24$ different
letters are fixed on the first step by no more than 
$(12g-6)(12g-6-2)(12g-6-4)$ possibilities. So, there are no
more than $(12g-6-24)$ posibilities for $a_2^{-1}b_2$,
no more than $(12g-6-26)$  for $b_2^{-1}c_2$ and
no more than $(12g-6-28$ for $c_2^{-1}b_2$ . There are no more
than $(12g-6-24)(12g-6-2-24)(12g-6-4-24)$ variants for fixing 24
letters on the second step. So, on the $k$-step we fix
$a_k^{-1}b_k, b_k^{-1}c_k$ and $c_k^{-1}a_k $ and their images.
There are no more than $(12g-6-24(k-1))(12g-8-24(k-1))(12g-10-24(k-1)$
posibilities for this. We have, that a genus $g$ Wicks form
can be constructed by no more than 
$$\prod_{m\in N, 12g-6-24(m-1)>0} (12g-6-24(k-1))(12g-8-24(k-1))
(12g-10-24(k-1))$$ ways. 
Lemma 3 follows now after some technical calcullations. 
\medbreak
\noindent
{\bf Lemma 4.} One may estimate the number of the Wicks forms with 
automorphisms of  order  $k$, $k\ge 7$ by the following expression
$$N_k(g)\le \left[\frac{12g-6}{7}\right]!.$$
\medbreak
\noindent
{\it Proof.} The proof is obvious.
\bigbreak
From the above Lemmas we obtain (\ref{sss}).

\hfill  $\Box$

Let us consider an element of the set of topological equivalence of the single-face nets, 
i. e. a word $U(g)\in W^r(g)$, $U(g)=a...b...c...$,
$\{ a, b, c, ... \}=A$ where $A$ is an alphabet. 
A corresponding local minimal net consists of arcs of geodesics 
$e_a, e_b,...$. 

Let $l_x$ be length of segment $e_x$ for  $x\in A$.
Then for any word $U(g)$ there is a system of numbers
$L_U=l_a,..., l_b,..., l_c,....$ where $l_x=l_{x^{-1}}$
and the subscripts put together the word $U(g)$.

Thus, the system of numbers defined above is a  description
of a local minimal net on a surface. On the other hand
there are different systems for any  net. So, we need to define 
now an equivalence of these systems. 
Two systems of numbers $L'_{U'}$ and $L_U$ are equivalent if
in some alphabet $U=U'$
and the corresponding components $l'_x$ and $ l_x$ are 
equal for any letter $x\in U$.  

What kind of  numbers in 
${\bf R}^{12g-6}_{+}$
correspond to a system of numbers $L_U$? 
In order to answer on this question 
we have to describe the set of $12g-6$-gons on the hypobolic plane
with the following conditions: all angles around vertecies 
are equal to $2\pi /3$,  sides of this polygon can be consequently 
labeled by 
letters of the word $U$ so that we can read this word after that 
 and for length $l_x$ of side $x$ holds $$l_x=l_{x^{-1}}, x\in A. 
\eqno (1)$$
So, there are $6g-3$ equalities of the type (1) because there are 
$6g-3$
different letters in  $U$. Let us now consider three 
different vertecies $A_1$, $A_2$, $A_3$ of such a $12g-6$-gon. We 
have then
three different nonintersect subwords $X_1$, $X_2$, $X_3$ of 
$U$. 
There are
analitic functions $\nu$, $\mu$ so that
$$\hat A_1=\nu (l_{X_{1}}, l_{X_2^{-1}}),$$      
$$\hat A_2=\nu (l_{X_{2}}, l_{X_3^{-1}}),$$      
$$\hat A_1=\nu (l_{X_1^{-1}}, l_{X_{3}}),$$      
$$|A_2A_3|=\mu (X_3),$$
$$|A_1A_2|=\mu (X_2),$$
$$|A_1A_3|=\mu (X_1).$$
It is well known that  there is  an analitic
function $\xi$ such that

$$|A_1A_2|=\xi (\hat A_3, \hat A_1, \hat A_2),  $$
$$|A_2A_3|=\xi (\hat A_1, \hat A_2, \hat A_3),  \eqno (2)$$
$$|A_3A_1|=\xi (\hat A_2, \hat A_3, \hat A_1).  $$

The equalities (1) and (2) give us a description of all
such polygons. The number of conditions is equal to $6g-3 + 3=6g$.

Thus, the set of local minimal nets $W_U$ corresponding to 
a cubic word $U$ by genus $g$ can be realized as a subset of 
${\bf R}^{12g-6}_{+}$ given by  equalities (1) and (2).

We obtain also the 
following theorem.
\bigbreak
\noindent
{\bf Theorem 1.2.} The set of all closed single-face local minimal nets  $W(g)$
 can be realized as a union of subsets $W_U$ in 
${\bf R}^{12g-6}_{+}$ given by  the conditions (1), (2)
where $U$ is an element of $W^r(g)$.
\bigbreak
\noindent
{\bf Commentar.} It is interesting to note that  the number 
of free parameters in the set $W_U$ is equal to $6g-6$, i. e. the  dimension
of the Teichm\"uller space.

\section
{\bf Exotic nets }

Here we suppose some examples of words corresponding to the regular $f$-face
regular nets, $f>1$. One can realize  such nets on the surfaces of
curvature $-1$. Let us consider $f$ regular $k$-gons
on the hypobolic plane with angles equal $2\pi /3$ 
and label the sides according these words.
Such  partitions of the hyperbolic plane with an action of the fundamental group
of $M^2_g$ according the following words  gives us  examples
of   surfaces with an "exotic" local minimal net.

Let $f$ be $2$ and $k$ be $12$. Then there is the following example
  
$$e^{-1}abca^{-1}db^{-1}c^{-1}d^{-1}elm^{-1},$$
$$mnfghf^{-1}ig^{-1}h^{-1}i^{-1}n^{-1}l^{-1}.$$

Let $f=3$  and $k=10$ :

$$a_1a_4a_5a_7a_1^{-1}a_8a_5^{-1}a_{14}^{-1}a_{10}^{-1}a_9^{-1}$$
$$a_{14}a_1^{-1}a_7^{-1}a_8^{-1}a_9a_{13}^{-1}a_{15}^{-1}a_{18}^{-1}
a_{19}^{-1}a_{11}^{-1}$$

$$a_{10}a_{11}a_{12}a_{15}a_{16}a_{18}a_{12}^{-1}a_{19}a_{16}^{-1}a_{13}$$

Let $f=4$  and $k=9$ :

$$a_1a_4a_5a_2a_3a_{15}^{-1}a_{18}^{-1}a_{19}^{-1}a_9^{-1}$$
$$a_2^{-1}a_6a_7a_1^{-1}a_8a_6^{-1}a_5^{-1}a_{11}^{-1}a_{10}$$
$$a_{11}a_4^{-1}a_7^{-1}a_8^{-1}a_9a_{12}a_{15}a_{16}a_{14}^{-1}$$

$$a_{14}a_{17}a_{18}a_{12}^{-1}a_{19}a_{17}^{-1}a_{16}^{-1}a_3^{-1}a_{10}$$

Let $f=6$  and $k=8$ :

$$a_1a_4a_5a_2a_3a_{20}a_{10}^{-1}a_9^{-1}$$
$$a_2^{-1}a_6a_7a_1^{-1}a_8a_6^{-1}a_5^{-1}a_{13}^{-1}$$
$$a_{13}a_4^{-1}a_7^{-1}a_8^{-1}a_9a_{14}a_{21}a_3^{-1}$$
$$a_{10}a_{11}a_{12}a_{15}a_{16}a_{23}^{-1}a_{22}^{-1}a_{14}^{-1}$$
$$a_{17}a_{18}a_{12}^{-1}a_{19}a_{17}^{-1}a_{16}^{-1}a_{24}a_{23}$$

$$a_{15}^{-1}a_{18}^{-1}a_{19}^{-1}a_{11}^{-1}a_{20}a_{21}a_{22}a_{24}^{-1}$$

Let $f=12$  and $k=7$ :

$$a_1a_2a_{20}a_{21}a_{22}a_{10}^{-1}a_9^{-1}$$
$$a_3a_4a_5a_{23}a_{24}a_{25}a_{20}^{-1}$$
$$a_6a_7a_3^{-1}a_2^{-1}a_{26}a_{27}a_{23}^{-1}$$
$$a_1^{-1}a_8a_6^{-1}a_5^{-1}a_{29}a_{28}^{-1}a_{26}^{-1}$$
$$a_4^{-1}a_{7}^{-1}a_8^{-1}a_9a_{30}a_{31}a_{29}^{-1}$$
$$a_{10}a_{11}a_{12}a_{13}a_{33}^{-1}a_{32}^{-1}a_{30}^{-1}$$
$$a_{14}a_{15}a_{16}a_{34}a_{35}a_{36}a_{33}$$
$$a_{17}a_{18}a_{14}^{-1}a_{13}^{-1}a_{37}a_{38}a_{34}^{-1}$$
$$a_{12}^{-1}a_{19}a_{17}^{-1}a_{16}^{-1}a_{40}a_{39}a_{37}^{-1}$$
$$a_{15}^{-1}a_{18}^{-1}a_{11}^{-1}a_{22}^{-1}a_{41}a_{40}^{-1}$$

$$a_{21}^{-1}a_{25}^{-1}a_{42}a_{35}^{-1}a_{38}^{-1}a_{39}^{-1}a_{41}^{-1}$$

$$a_{42}^{-1}a_{24}^{-1}a_{27}^{-1}a_{28}a_{31}^{-1}a_{32}a_{36}^{-1}$$

\end{document}